\documentclass[12pt]{amsart}
\usepackage{amsmath}
\usepackage{amssymb}
\usepackage{anysize}
\marginsize{3cm}{3cm}{3cm}{3cm} 

\def\appendix{\par}  

\newtheorem{claim}{Claim}
\theoremstyle{plain}
\newtheorem{theor}{Theorem} 
\newtheorem{lemma}[theor]{Lemma} 
\newtheorem{propo}[theor]{Proposition} 
\newtheorem{corol}[theor]{Corollary}  
\theoremstyle{definition} 
\newtheorem{defin}[theor]{Definition}  
\newtheorem{quest}[theor]{Question} 

\newcommand{\eq}{{{\mathfrak e \mathfrak q}}}

\newcommand{\body}[1]{[#1]}
\newcommand{\finseq}{{}^{\stackrel{\omega}{\smile}}\omega}
\newcommand{\binseq}{{}^{\stackrel{\omega}{\smile}} 2}

\DeclareMathOperator{\cov}{\mathbf{cov}}
\newcommand{\AND}{\text{ and }}

\newcommand{\OR}{\text{ or }} 
\newcommand{\card}[1]{\left\lvert #1 \right\rvert}

\newcommand{\forces}[1]{\Vdash \mbox{``} #1 \mbox{''}}

\newcommand{\SetOf}[2]{\left\{#1 \ \left| \ #2 \right.\right\}}

\newcommand{\Reals}{{\mathbb R}}

\newcommand{\Poset}{{\mathbb P}}

\newcommand{\Integers}{{\mathbb Z}}
\newcommand{\Naturals}{{\mathbb N}}
\newcommand{\ba}{{\mathbb B}}

\newcommand{\presup}[2]{\, ^{#1} \! #2} 
\newcommand{\wpresup}[1]{\presup{\stackrel{\omega}{\smile}}{#1}} 

\newcommand{\wfomom}{\wpresup{\omega}}

\title[Translates of a closed nowhere dense set]{The number of translates 
of a closed nowhere dense set required to cover a Polish group}
\author[A.~W.~Miller]{Arnold~W.~Miller}
\address{Department of Mathematics, University of Wisconsin,
Van~Vleck Hall, 480 Lincoln Drive,
Madison, Wisconsin, USA 53706-1388}
\curraddr{}
\email{miller@math.wisc.edu \hfill www.math.wisc.edu/$\sim$miller}
\author[J.~Stepr\={a}ns]{Juris Stepr\={a}ns}
\address{Department of Mathematics, York University,
4700 Keele Street,
Toronto, Ontario, Canada M3J 1P3}
\curraddr{}
\email{steprans@yorku.ca \hfill www.math.yorku.ca/$\sim$steprans}
\thanks{Research for this paper was partially supported by NSERC of
Canada.}
\keywords{translations, nowhere dense sets, Polish groups,
cardinal invariants}
\subjclass{03E17}

\begin{document}
\begin{abstract}
For a Polish group ${\mathbb{G}}$ let $\cov_{\mathbb{G}}$ be the minimal
number of translates of a fixed closed nowhere dense
subset of ${\mathbb{G}}$ required to cover ${\mathbb{G}}$.  For  
many locally compact ${\mathbb{G}}$ this cardinal is known to	 
be consistently larger than $\cov({\mathcal M})$ which is the smallest
cardinality of a covering of the real line by meagre sets.  It is shown that
for several
non-locally compact groups   $\cov_{\mathbb{G}}=\cov({\mathcal M})$. For
example the equality holds for the group of permutations of the
integers, the additive group of a  
separable Banach space with an unconditional basis and the group of
homeomorphisms of various compact spaces. 
\end{abstract}
\maketitle
\bibliographystyle{plain}

The notion of translation invariants corresponding to the usual
invariants of the continuum has been considered by various researchers
and an introductory survey can be found in \S2.7 of
the monograph \cite{MR96k:03002} by Bartoszy\'nski and Judah.  The key
definition for the purposes of this 
article is the cardinal they denoted by $\cov^*({\mathcal M})$.
It is the least cardinal of a set $X\subseteq\Reals$ such 
that there is some meagre set $M\subseteq \Reals$ such
that $X+M = \Reals$. It is asserted that the value of
$\cov^*({\mathcal M})$ will be the same if the group $(\Reals, +)$ is
replaced in this definition by the Cantor set with its natural Boolean
operation or an infinite product
of finite cyclic groups.  The goal of this note is to initiate
a study of translation invariants for arbitrary Polish groups by
establishing that not all Polish groups yield the same invariants and
posing various questions which arise from this observation. Since many
of the interesting questions in this area concern non-locally compact
groups the measure-theoretic version is not easily formulated
and, therefore, only the topological version will be considered.
 
  Throughout, the statement that ${\mathbb G}$ is a group will mean
  that ${\mathbb G} = (G,\cdot,{\/}^{-1})$ but for $x$ and $y$ in $G$
  the operation $x\cdot  y$ will usually be abbreviated to
  $xy$. Similarly, if $A\subseteq G$ and $B\subseteq G$ then $AB$ will
  denote the set $\SetOf{xy}{x\in A \AND y\in B}$. If $A = \{a\}$ then
  $\{a\}B $ will be abbreviated to $aB$.

To begin, a generalization of $\cov^*({\mathcal M})$ will be defined
for arbitrary group actions. 
\begin{defin}
Let $\mathbb{G}$ be a group acting on a Polish space $X$
with the action denoted by $\alpha:G \times X \to X$.
Define $\cov_{\alpha} $ to be the least cardinal of a set
$Q\subseteq G$ such that there is some closed nowhere  dense set
$C\subseteq X$ such that $\alpha(Q,C) = \SetOf{\alpha(h,c)}{h\in Q
  \AND c\in C} = X$.
Define $\cov^*_{\alpha} $ to be the least cardinal of a set
$Q\subseteq G$ such that there is some meagre set
$M\subseteq X$ such that $\alpha(Q,M) = X$.
 In the special case of a Polish group
$\mathbb{G}$ acting on 
itself by left translation these invariants will  be denoted by
$\cov_{\mathbb{G}}$ and $\cov^*_{\mathbb{G}}$.  
Note that  the definitions of $\cov_{\mathbb{G}}$ and $\cov^*_{\mathbb{G}}$
do not change if left translation is replaced by right translation since
$(Q C)^{-1}=C^{-1} Q^{-1}$. 
\end{defin}

The first observation applies to arbitrary groups.

\begin{propo}\label{p:m=c}\label{action}
If $\mathbb{G}$
is an arbitrary group, $X$ is $\sigma$-compact, 
second countable and has no isolated points, 
$\alpha: G \times X \to X$ is a group action all of whose orbits are 
dense in $X$ and $\alpha(g,\cdot)$ is continuous for each fixed $g\in G$, then
$\cov_{\alpha}  = \cov^*_{\alpha}$.
\end{propo}
\begin{proof}
It suffices to show that if $M$ is meagre then
there is a closed nowhere dense set $C$ such that 
$\alpha(Y,C)\supseteq M$ for some
countable set $Y$.  Since $X$ is $\sigma$-compact it may be assumed that
$M= \bigcup_{n=0}^\infty  D_n$ where each $D_n$ is compact
and nowhere dense.  Let $\{B_n\}_{n=0}^\infty$
be a base for $X$ consisting of non-empty open sets.  
Select $Y_k$ and $U_k$ by induction on $k$ such that 
\begin{itemize}
  \item $Y_k$ is a finite subset of $G$ 
  \item $U_k\subseteq B_k$ is a non-empty open set
  \item if $j \leq k$ and $d\in D_j$ then there is some $ y\in Y_j$
such that $   \alpha(y, d)\notin\bigcup_{i=0}^k U_i$
  \item the closure of $\bigcup_{j=0}^k U_j$ is not equal to $X$
\end{itemize}
If this can be done then let $C$ be the complement of
$\bigcup_{i=0}^\infty U_i$ and take $Y$ countable including
$\bigcup_{n=0}^\infty Y_i$ and closed under inverses.
Then given $x\in M$ there is some $n$ such that $x\in D_n$ and hence
there exists $y\in Y_n$ such that $\alpha(y, x)\in C$ and so 
$$x =
\alpha(y^{-1},\alpha(y, x)) \in \alpha(\{y^{-1}\},C)\subseteq \alpha(Y,C).$$

To carry out the induction suppose this has been done for $i<k$.
Since $\alpha$ is a continuous group action, $\alpha(\{g\},D_j)$ is
closed nowhere dense for any $g\in G$ and so it is possible to choose
a non-empty open $U_k\subseteq B_k$ disjoint from
$$\bigcup_{i < k}\bigcup_{y\in Y_i}\alpha(\{y\}, D_i).$$
Moreover, it is possible to choose $U_k$ so that there is some non-empty
open set $W$ witnessing the last induction clause; in  other words,
$W\cap \bigcup_{i=0}^kU_i = \emptyset$. Since $D_k$ is compact and
$\alpha(\{d\},G)$ is dense for each $d\in D_k$ there is a
finite set $Y_k\subseteq G$ such that for each $d\in D_k$ there is
some $y \in Y_k$ such that $\alpha(y,d)\in W$.
\end{proof}

\begin{propo}\label{covpolish}\label{cov}
$\cov_{\mathbb{G}}  = \cov^*_{\mathbb{G}}$ for every non-discrete Polish 
$\mathbb{G}$.
\end{propo}

\begin{proof}
  If $\mathbb{G}$ is compact, then apply Proposition~\ref{action}. On
  the other hand, Birkhoff and Kakutani proved that every Polish group
  admits a left invariant\footnote{In other words,
    $d(y,z)=d(x y,x z)$ for every $x,y,z\in G$.  For a proof see
    \cite{MR2000k:03097} \S7.1 page 115.}  metric
 $d$ inducing the same topology.
  While this metric $d$ might not be complete, it is shown in
  Corollary~1.2.2 of \cite{MR98d:54068} that the metric $D$ defined by
  $D(x,y)= d(x,y)+ d(x^{-1},y^{-1})$ is a complete metric compatible  
  with the topology.  However $D$ might not be left or right invariant.
  Nevertheless, since $\mathbb{G}$ is not compact and $D$ is complete
  it follows that $D$ is not totally bounded;  in other words, there
  is some $\epsilon > 0$ and a sequence $\{x_n\}_{n\in\omega}$ such
  that $D(x_n,x_m) > \epsilon$ for all pairs of distinct integers $n$
  and $m$.  A routine application of Ramsey's Theorem yields an
  infinite subset $S\subseteq \omega$ such that either $d(x_n,x_m) >
  \epsilon/2$ for all pairs of distinct integers $n$ and $m$ in $S$ or
  $d(x_n^{-1},x_m^{-1}) > \epsilon/2$ for all pairs of distinct
  integers $n$ and $m$ in $S$.  From the left invariance of the metric
  $d$ it follows that letting $U$ be the $\epsilon/3$ ball around the
  identity either $\{x_n U\}_{n\in\omega}$ are pairwise disjoint sets
  or
$\{x_n^{-1} U\}_{n\in\omega}$ are pairwise disjoint sets. Without
  loss of generality assume the former alternative.

As in the argument for Proposition~\ref{action}, it suffices to prove the 
following:
For any family $\{C_n\}_{n\in\omega}$ of nowhere
dense subsets of ${\mathbb{G}}$ there exists a countable 
$Q\subseteq {\mathbb{G}}$ and a nowhere dense $C$ such that
$\bigcup_{n\in\omega} C_n\subseteq Q C$.
Choose $\{q_n\}_{n\in\omega}\subseteq{\mathbb{G}}$  
so that $\bigcup_{n\in\omega}q_n U={\mathbb{G}}$. 
This is possible since $\SetOf{x U}{x\in G}$ covers $G$ and
$G$ is Lindelof.  Let $k:\omega\times\omega\to\omega$
be a bijection and define $r_{n,m}$ so that $r_{n,m}
q_n=x_{k(n,m)}$ and define
$$C=\bigcup_{m\in\omega}\bigcup_{n\in\omega}
 r_{n,m}(q_n U\cap C_m).$$
Note that $r_{n,m} (q_n U\cap C_m)\subseteq x_{k(n,m)} U$ 
and since $\{x_n U\}_{n\in\omega}$ are disjoint, $C$ is nowhere dense. 
On the other hand
if $Q=\{r_{n,m}^{-1}\}_{n,m\in\omega}$, then $Q C$ contains 
$q_n U\cap C_m$ for each $n$ and $m$ and hence $C_m\subseteq Q C$. 
\end{proof}

The proof of the following fact can be found in
\cite{MR96k:03002} in Lemma~2.4.2.  

\begin{theor}\label{t.2.4.2}
$$\cov({\mathcal M})=\min\SetOf{\card{\mathcal F}}{
{\mathcal F}\subseteq \omega^\omega\AND\forall g\in \omega^\omega\ 
\exists f\in{\mathcal F}\ \forall n\ f(n)\neq g(n)}.$$
\end{theor}

\begin{corol}
 $\cov_{\Integers^\omega}=\cov({\mathcal M})$.
\end{corol}
\begin{proof}
For any non-discrete Polish group ${\mathbb{G}}$ the inequality
$\cov_{\mathbb{G}}\geq\cov({\mathcal M})$ holds.
Let 
$$C=\SetOf{f\in \Integers^\omega}{\forall n\; f(n)\neq 0}.$$
Note that $C$ is closed nowhere dense.  
Take ${\mathcal F}\subseteq\omega^\omega\subseteq \Integers^\omega$ 
such that  $\card{\mathcal F}=\cov({\mathcal M})$ and such that 
$$\forall g\in \Integers^\omega\;\exists f\in {\mathcal F}\;\forall n\;
f(n)\not=g(n).$$
But this means that ${\mathcal F}+C=\Integers^\omega$. Hence
$\cov_{\Integers^\omega}\leq \card{\mathcal F}=\cov({\mathcal M})$.
\end{proof}

To generalize this to other groups the following lemma will play a key role.
It is stated using standard notation about trees.
If $T\subseteq \wfomom$ is a tree then
$$\body{T} = \SetOf{f}{\forall n\in \omega\ f\restriction n\in T}$$ 
and $T$ is  infinite branching if for each $t\in T$ the set
$\Sigma(t) = \SetOf{n\in\omega}{t^\frown (n)\in T}$
 is infinite.
\begin{lemma}\label{l.2.4.2}
If $T\subseteq\finseq$ is an infinite
branching tree then there is ${\mathcal F}$ such
that 
\begin{itemize}
\item $\card{\mathcal F} = \cov({\mathcal M})$
\item ${\mathcal F}\subseteq \body{T}$ 
\item   for each $g:\omega\to \omega$ there is $f\in {\mathcal F}$
  such that $f(n) \neq g(n)$ for all integers $n$.
\end{itemize}
\end{lemma}
\begin{proof}
Using Theorem~\ref{t.2.4.2} let ${\mathcal F}^*\subseteq \omega^\omega$
be a family of cardinality $\cov({\mathcal M})$                    
such that for every $g\in \omega^\omega$ there is $f\in {\mathcal  
F}^*$ such that $f(n)\neq g(n)$ for every integer $n$.             
By choosing an infinite branching subtree of $T$ if necessary, it may be  
assumed that $\Sigma(t)\cap \Sigma(s) = \emptyset$ unless $s=t$.
Let $\psi_t:\omega\to \Sigma(t)$ be a bijection for each $t\in T$. For
any $f:\omega\to\omega$ define $f_\psi$ by inductively 
setting $$f_\psi(n) = \psi_{f_\psi\restriction n}(f(n))$$ and 
note that $f_\psi\restriction n \in T$ for all $n$.
Then let
${\mathcal F} = \SetOf{f_\psi}{f\in {\mathcal F}^*}$.

Let $\Sigma=\cup_{s\in T}\Sigma(s)$.
Given $g:\omega \to \omega$ and $k\in \omega$ if $g(k)\in\Sigma$
let $s_k \in T$ be the unique node in $T$ with $g(k)\in \Sigma(s_k)$.

Now define $G:\omega\to\omega$ by
$$G(k) = 
\begin{cases}
  \psi^{-1}_{s_k}(g(k)) & \text{if $g(k)\in\Sigma$ }\\
   0 &\text{otherwise }
\end{cases}$$
and choose $f\in{\mathcal F}^*$ such that $f(n) \neq G(n)$ for all
$n$. It suffices to show that $f_\psi(n)\neq g(n)$ for all
$n$.  Assume that $f_\psi(n) = g(n)$. Then
necessarily $g(n)=f_\psi(n)\in\Sigma(f_\psi\restriction n)$ 
and so $s_n=f_\psi\restriction n$ and therefore
$$G(n)=\psi^{-1}_{s_n}(g(n))=\psi^{-1}_{f_\psi\restriction
  n}(f_\psi(n))=f(n)$$  
which is a contradiction.

\end{proof}

\begin{theor}\label{m:tn}
Let ${\mathbb G}$ be a Polish group such that there are
$\{B_k, A_k^j\}_{j,k\in\omega}$ such that:
\begin{enumerate}
\item $B_k$ and $A_k^j$ are all subsets of $G$ 
\item there is an infinite branching tree $T\subseteq
  {}^{\stackrel{\omega}{\smile}}\omega$ such that
  $\bigcap_{k\in\omega}A^{b(k)}_k\neq \emptyset$ for $b\in\body{T}$
  \item $\bigcup_{k\in\omega}B_k$ is dense open
  \item $(A_k^i B_k )\cap (A_k^j B_k) = \emptyset$ unless $i=j$
\end{enumerate}
then $\cov_{\mathbb G} = \cov({\mathcal M})$.
\end{theor}

\begin{proof}
Using Lemma~\ref{l.2.4.2} let ${\mathcal F}$ be a set of branches
through $T$ such that 
 $\card{\mathcal F} = \cov({\mathcal M})$
and   for each $g:\omega\to \omega$ there is $f\in {\mathcal F}$
  such that $f(n) \neq g(n)$ for all integers $n$.
  For $b\in\body{T}$ let $b^* \in \bigcap_kA_k^{b(k)}$ and
  let $X= \SetOf{b^*}{b\in{\mathcal F}}$.
Let $C= G\setminus \bigcup_{k}B_k$. Then $C$ is closed and nowhere
dense and $\card{X} = \cov({\mathcal M})$. Hence it suffices to
show that $X C = G$.

To this end let $g\in G$ and define $\Gamma:\omega \to \omega$ such that
$g\in A_n^{\Gamma(n)} B_n$ if there is any $i$
with  $g\in A_n^i B_n$. By
Hypothesis~3 the choice of $\Gamma(n)$ is unique if it exists at all.
Let $f\in {\mathcal F}$ be such that $f(n) \neq \Gamma(n)$ for all integers
$n$. If $g\notin X C$ then $g\notin f^* C$  and so
there is some integer $m$ such that $g\in f^* B_m\subseteq
A_m^{f(m)} B_m$. Hence $\Gamma(m) = f(m)$ which is impossible.
\end{proof}

For example, for the group $\Integers^\omega$ one could
take $A^j_k=\SetOf{\alpha}{\alpha(k)=j}$, $B_k=A_k^0$, and
$T=\finseq$ or, for the group $\Reals^\omega$ take 
$B_k=\SetOf{\alpha\in \Reals^\omega}{|\alpha(k)|<1/3}$ and
let $A^j_k$ and $T$ be the same.

\begin{corol}
  If ${\mathbb G}$ is the group of all permutations of the integers
  then $\cov_{\mathbb G} = \cov({\mathcal M})$.
\end{corol}

\begin{proof}
Let $A_i^k = \SetOf{p}{p(2i) = 2k}$ and $B_k=A_k^k$. Letting
$T$ be the tree of all one-to-one sequences satisfies the hypotheses of 
Theorem~\ref{m:tn}.  
\end{proof}
For any compact metric space $X$ let ${\mathbb H}(X)$ be the group of
autohomeomorphisms of $X$ using composition as the group operation and the
topology induced by the uniform metric.

\begin{corol}\label{c:hui}
$\cov_{{\mathbb H}([0,1])} = \cov({\mathcal M})$.
\end{corol}
\begin{proof}
Construct a family of open intervals 
$\SetOf{I_s}{s\in\Integers^{<\omega}}$ so that: 
\begin{itemize}
 \item $I_{\emptyset}=(0,1)$
 \item $I_s =\bigcup_{n\in\Integers}\overline{I_{s^\frown (n)}}$ and the
 right-hand endpoint of $I_{s^\frown(n)}$ is the left-hand endpoint of
 $I_{s^\frown(n+1)}$
 \item the length of $I_s$ is less than $\frac{1}{|s|}$.
\end{itemize}
In other words,  $I_s$ is partitioned into contiguous subintervals in order
type $\Integers$ and their endpoints. 
Define for $i,k\in \omega$ the open sets
$$U_k^i=\cup\SetOf{I_s }{ s\in\Integers^{k+1} \AND s(k)=i}$$

Fix any $p_0\in (0,1)$. Define $$B_k=\SetOf{h\in{\mathbb H}}{
  h(p_0)\in U_k^0}$$ 
and let $U=\bigcup_{k\in\omega}U_k^0$ and note that $U$ is open dense
and 
$$\bigcup_{k\in\omega}B_k=\SetOf{h\in{\mathbb H}}{ h(p_0)\in U}$$
is open dense in $\mathbb H$.   
For any $k,i\in\omega$ define
$$A_k^i=\SetOf{h\in{\mathbb H}}{ h(U_k^0)\subseteq U_k^i} .$$ 
It is clear that $\{A_k^i\circ B_k\}_{i\in\omega}$ are pairwise 
disjoint because if $h\in A_k^i\circ B_k$ then
$h(p_0)\in U_k^i$ and $\{U_k^i\}_{i\in\omega}$ are pairwise disjoint.

Now let $T=\finseq$ and suppose $b\in\body{T}$. 
In order to find $h\in \bigcap_{k\in\omega}A_k^{b(k)}$ 
let $a_s$ denote the left-hand endpoint of $I_s$ and define
$Q=\SetOf{a_s}{s\in\Integers^{<\omega}}$.  Note that $Q$ is
dense in $[0,1]$ and its ordering 
is exactly the same as the lexicographical ordering $<_{lex}$
on $\Integers^{<\omega}$:
Define $s <_{lex} t$ if and only if
$s\subset t$ or there exists $i$ such that $s(i)<t(i)$ and 
$s\restriction i=t\restriction i$.
Then $a_s<a_t$ if and only if $s<_{lex} t$. 
for all $s,t\in \Integers^{<\omega}$.  

For each $s\in \Integers^n$ define $s+b=t\in \Integers^n$ by
pointwise addition and define
$h(a_s)=a_{s+b}$.  Clearly the mapping
$s\mapsto s+b$ is a bijection preserving the
lexicographical order and so 
$h:Q\to Q$ is an order preserving bijection.
Thus it extends uniquely to an order preserving  
bijection $h^*$ on $[0,1]$.  Note that 
$h(I_{s^\frown(0)})\subseteq I_{t^\frown (b(k))}$ for any $t$, $s$ and
$b$ where 
$t=s+b$ and $k=|s|$.   So
$h^*(U_k^0)\subseteq U_k^{b(k)}$ for each $k$. Thus
$h^*\in \bigcap_{k\in \omega}A_k^{b(k)}$.
\end{proof}
\begin{corol}\label{c:g}
Suppose there exists $\pi:X\to [0,1]$ which is
\begin{enumerate}
 \item continuous, onto, and open
 \item $\forall h\in{\mathbb H}([0,1])\;\;\exists \hat{h}\in{\mathbb H}(X)\;\;
  \forall x\in X\;\; \pi(\hat{h}(x))=h(\pi(x))$
 \item if $U\subseteq (0,1)$ is open dense, then there exists $p_0\in X$
 such that $$\SetOf{h\in{\mathbb H}(X)}{ h(p_0)\in \pi^{-1}(U)}$$ is open dense
\end{enumerate}
then $\cov_{{\mathbb H}(X)}=\cov({\mathcal M})$.\end{corol}
\begin{proof}
  Using the same set \/ $U = \bigcup_{k\in\omega}U_k^0 \subseteq
 (0,1)$ as defined 
 in the proof of Corollary~\ref{c:hui} let $p_0\in X$ be
 such that 
$$\SetOf{h\in{\mathbb H}(X)}{ h(p_0)\in \pi^{-1}(U)}$$ is open dense.
 Define $$B_k=\SetOf{h\in{\mathbb H}(X)}{
 \pi\circ  h(p_0)\in U_k^0}$$ 
and note that 
$$\bigcup_{k\in\omega}B_k=\SetOf{h\in{\mathbb H}(X)}{ h(p_0)\in \pi^{-1}(U)}$$
is open dense in ${\mathbb H}(X)$.   
For any integers $k$ and $i$ define
$$A_k^i(X)=\SetOf{h\in{\mathbb H}(X)}{ \pi\circ h\left(
  \pi^{-1}(U_k^0)\right)\subseteq U_k^i}$$   
reserving the notation $A_k^i$ for the sets defined in the proof of
Corollary~\ref{c:hui}. 
It is clear that $\{A_k^i(X)\circ B_k\}_{i\in\omega}$ are pairwise 
disjoint because if $h\in A_k^i(X)\circ B_k$ then
$\pi\circ h(p_0)\in U_k^i$ and $\{U_k^i\}_{i\in\omega}$ are
  pairwise disjoint. 

Now let $T=\finseq$ and suppose $b\in\body{T}$. 
It has already been established in the proof of Corollary~\ref{c:hui}
that there is $h\in \bigcap_{k\in\omega}A_k^{b(k)}$. 
Then it is readily verified that $\hat{h} \in
\bigcap_{k\in\omega}A_k^{b(k)}(X)$.  
\end{proof}

\begin{corol}
  If $X$ is either $[0,1]^k$ for $k\in \omega$ or  the infinite
  dimensional Hilbert cube 
  $[0,1]^\omega$ then  
$\cov_{{\mathbb H}(X)} = \cov({\mathcal M})$.
\end{corol}
\begin{proof}
Use Corollary~\ref{c:g} applied to the mapping which projects 
a sequence to its first coordinate. Given $h\in{\mathbb H}([0,1])$
take $\hat{h}\in {\mathbb H}([0,1]^k)$ to be defined
by 
$$\hat{h}(x_0,x_1,x_2,\ldots)=(h(x_0),x_1,x_2,\ldots).$$ 
Let $p_0=(1/2,1/2,\ldots,1/2)$ and note that 
by a classical result of L.E.J. Brouwer (invariance of
domain) 
$h(p_0)$ is in the interior of $[0,1]^k$ and so its first coordinate
is neither $0$ or $1$.  It is easy to see that condition (3) of
Corollary \ref{c:g} holds.  In the case of $[0,1]^\omega$ it is possible
that the first coordinate of $h(p_0)$ is $0$ or $1$ since the Hilbert
cube is homogeneous.  However, one of the standard proofs of homogeneity,
see Lemma 6.1.4 page 252 in Van Mill \cite{vanmill}, 
shows that given any
$(u_n:n\in\omega)\in [0,1]^\omega$,
there exists a small $k\in\hom([0,1]^{\omega})$
with $k(u_n:n\in\omega)=(v_n:n\in\omega)$
and $v_0\in (0,1)$.
\end{proof}
Since the homeomorphism group of the  Hilbert cube is universal --- in the
sense that any Polish group is homeomorphic to one of its closed
subgroups --- the following question is of central importance in this area. 
\begin{quest}
  If ${\mathbb G}$ is a non-$\sigma$-compact, closed subgroup of 
${\mathbb H}([0,1]^\omega)$ does the equality
  $\cov_{\mathbb G} = \cov({\mathcal M})$ hold?
\end{quest}

\begin{corol}
If $S^n$ is the $n$ dimensional Euclidean sphere then 
$$\cov_{{\mathbb H}(S^n)} = \cov({\mathcal M}).$$ Similarly
for any compact metric space $X$ we have that 
$$\cov_{{\mathbb H}(S^n\times X)} = \cov({\mathcal M}).$$
\end{corol}
\begin{proof}
Use Corollary~\ref{c:g} applied to the mapping which projects the
sphere onto one of its diameters.  
\end{proof}

\begin{corol}\label{banach}
  If ${\mathbb G}$ is the additive group of a separable Banach space
  with an unconditional basis
  then $\cov_{\mathbb G} = \cov({\mathcal M})$.
\end{corol}
\begin{proof}
Let $\{e_i\}_{i\in\omega}$ be an unconditional basis for the Banach
space $\ba$. Recall that  this implies that
for each $S\subseteq \omega$ the projection map
$$P_S\left(\sum_{n<\omega}\alpha_ne_n\right)=\sum_{n\in
  S}\alpha_ne_n$$
is a well-defined continuous linear operator.

Let $\{S_j\}_{j\in\omega}$ be pairwise disjoint infinite subsets of $\omega$.
Choose $\{\delta_i\}_{i\in\omega}$ positive reals such that
$\sum_{i=0}^\infty\delta_i < \infty$. 
For each $k$ the range of $P_{S_k}$ is an infinite dimensional Banach
space and hence  the ball of diameter
$\delta_k$ around the origin is not compact and therefore not totally
bounded. In other words, it is possible to find some $\epsilon_k > 0$ and
$\{u_k^i\}_{i\in\omega}$ contained in the range of $P_{S_k}$ such that $\|u_k^i\| < \delta_k$ and
$\|u_k^i-u_k^j\| > \epsilon_k$ 
for distinct $i$ and $j$.
Let
$$B_k=\SetOf{b\in \ba}{\|P_{S_k}(b)\|<\epsilon_k/2}
\AND A_k^{j} =\SetOf{b\in \ba}{P_{S_k}(b) = u_k^j}.$$
Finally, let
$T=\bigcup_{n<\omega}\prod_{k<n}S_k$ and note that $T$ is infinite branching.
 
To check that the other hypotheses of Theorem~\ref{m:tn} are satisfied
note that  
$B_k$ is open because $P_{S_k}$ is continuous. To see that
$\bigcup_kB_k$ is dense 
note that the finite linear sums of the basis form a dense subset and if
$b=\sum_{i=0}^N\alpha_ie_i$ and
$S_k$ is disjoint from $\{0,1,\ldots,N\}$ then $P_{S_k}(b) = {\bf 0}$.
To see that 
$\bigcap_k A^{b(k)}_k$ is non-empty for $b\in \body{T}$ 
note that $b^*=\sum_{k=0}^\infty  u_k^{b(k)}$ is a convergent series
by the choice of the
$\delta_i$.
The sets $\{A_k^j + B_k\}_{j\in \omega}$ 
are pairwise disjoint by the choice of $\epsilon_k$
\end{proof}

\begin{quest}
Is Corollary~\ref{banach} 
true for Banach spaces without an unconditional basic sequence?
\end{quest}

\begin{defin}
For $f:\Naturals \to \Naturals$ define $\eq(f)$ to be the least      
cardinal of a set $X\subseteq \prod_{i=0}^\infty f(i)$ such that for    
all $g\in\prod_{i=0}^\infty f(i)$ there is $x\in X$ such that           
$x(i)\neq g(i)$ for every $i\in \Naturals$. Let $\eq$ be the    
minimum of all $\eq(f)$ such that $\lim_{i\to\infty}f(i) = \infty$.  
\end{defin}

\begin{propo}\label{torus}
(1) $\cov_{\mathbb G} = \eq$ when ${\mathbb G}$ 
is either the product of finite cyclic
groups or $\Reals$.  
\par (2) Let ${\mathbb G}=\mathbb{T}=\Reals/\Integers$ be the circle group or
 any finite dimensional torus ${\mathbb{T}^n}$.  Or
 let ${\mathbb G}$ product of countably many nontrivial
 finite groups. Then $\cov_{\mathbb G} = \eq$.
\end{propo}

\begin{proof}
Assertion~(1) is proved in \cite{MR96k:03002} and their argument
generalizes to (2). 
\end{proof}

\begin{lemma}\label{inequal}
Suppose $\mathbb G$ and $\mathbb H$ are topological groups and 
$h:\mathbb G \to \mathbb H$ is a   
homomorphism onto $H$. 
\par (a) If $h$ is open and continuous, then                                          
$\cov_{\mathbb H} \geq \cov_{\mathbb G}$.  
\par (b) If $h$ takes
meagre sets to meagre sets, then 
$\cov^*_{\mathbb H} \leq \cov^*_{\mathbb G}$.
\end{lemma}
\begin{proof}
  
  (a) Since $h$ is open and continuous, the pre-image of a nowhere
  dense set is nowhere dense.  Let $C\subseteq \mathbb H$ be nowhere
  dense and $Y\subseteq \mathbb H$ have the property that $|Y|=
  \cov_{\mathbb H}$ and $YC=H$.  Then $B=h^{-1}(C)$ is nowhere dense
  in $\mathbb G$ and choose $Y$ so that $h(Y)=X$ and $|X|=|Y|$.  Since
  $h$ is homomorphism it follows that $h(YB)=XC=H$.  In order to see
  that $YB=G$ note that for any $g\in G$ there exists $y\in Y$ and
  $b\in B$ with $h(g)=h(yb)$.  It follows that $h(y^{-1}g)=h(b)\in C$
  so $y^{-1}g\in B$ and so $g=yb^\prime$ for some $b^\prime\in B$.

(b) Suppose $h$ takes meagre sets to meagre sets and let $M\subseteq \mathbb G$
be meagre such that $\cov^*_{\mathbb G}=|X|$ for some $X$ with $XM=G$.  Then
$h(X)h(M)=H$ and so  $\cov^*_{\mathbb H} \leq |X|=\cov^*_{\mathbb G}$

\end{proof}

So, for example, if $\mathbb G$ and $\mathbb H$ are Polish groups then  
$\cov_{\mathbb H\times \mathbb G} \leq \cov_{\mathbb G}$.  

\begin{propo}
For any integer $n$ the equality $\cov_{\Reals^n} = \eq$ holds.
\end{propo}
\begin{proof}
The quotient mapping from $\Reals^n $ to $\Reals^n/\Integers^n$ is
both open and it takes nowhere dense
sets to meagre sets.  It follows from Lemma \ref{inequal} that
$\cov_{\Reals^n}\leq\cov_{\Reals^n/\Integers^n}$ and
$\cov^*_{\Reals^n}\geq\cov^*_{\Reals^n/\Integers^n}$.
  
Since
$\Reals^n/\Integers^n$ is isomorphic to $\mathbb{T}^n$ it follows from
Proposition~\ref{torus} and Corollary~\ref{cov}
that $\cov_{\Reals^n} = \eq$.
\end{proof}

There are many models of set theory where $\cov({\mathcal M})<\eq$;
examples can be found using modified Silver forcing as in \S7 of 
\cite{MR84e:03058a}.
\nocite{MR84e:03058b}

\begin{theor}\label{consistent}
It is relatively consistent with ZFC that 
\begin{itemize}
\item $2^{\aleph_0} = \aleph_2$,
\item $\cov_{\mathbb
    G}=\aleph_2$ for every infinite compact group ${\mathbb G}$
\item and $\cov({\mathcal M})=\aleph_1$.
\end{itemize}
 In fact, the last equality will from the fact 
 that ${\mathfrak d}=\aleph_1$ holds in the model constructed. 
\end{theor}
\begin{proof}
Fix ${\mathbb G}$ an infinite compact group and
$C\subseteq G$ closed nowhere dense. Let $h:2^\omega\to G$ be a 
continuous onto map. 
Recall the notation preceding Lemma~\ref{l.2.4.2} concerning trees.
If $p$ is a tree and $s\in p$ define 
 $p_s=\SetOf{t\in p}{ t\subseteq s \OR s\subseteq t}$.
\begin{defin}
Define the partially ordered
set $\Poset=\Poset_{(h,C,{\mathbb G})}$ as follows:  
$p\in\Poset$ if and only if
$p\subseteq \binseq$ is a perfect tree such that
$xC$ is relatively nowhere dense
in $h(\body{p_s})$ for every $x\in G$ and $s\in p$. 
This is equivalent to saying 
for every $x\in G$ and $s\in p$ there exists $t\in p$ with
$s\subseteq t$ and $h(\body{p_t})\cap xC=\emptyset$.
The ordering on $\Poset$ is inclusion considering a subtree as a stronger
condition.
\end{defin}
First note that the forcing is not vacuous.
For $s\in \binseq$ let $B_s=\SetOf{x\in 2^\omega}{ s\subseteq x}$.
If $p=\SetOf{s\in \binseq}{ h(B_s) \text{ is non-meagre in } G}$ then 
$p\in \Poset$.  
It will be shown that the forcing has the following three properties:
\begin{itemize}
\item[(a)] $\Poset_{(h,C,{\mathbb G})}$ is a proper forcing notion.  
\item[(b)] Forcing with $\Poset_{(h,C,{\mathbb G})}$ adds no $\leq^*$
 unbounded  sequence in  $\omega^\omega$.
\item[(c)] In the generic extension there exists $z\in G$ such that $z\notin XC$ where
$X$ is the set of elements of $G$ in the ground model.     
\end{itemize}

\begin{claim}
Given any $p\in \Poset$ and $n<\omega$ there exists $N$ with $n<N<\omega$ 
such that for every $s\in 2^n\cap p$ and $x\in G$ there exists 
$t\in 2^N\cap p_s$ such that $h(\body{p_t})\cap xC=\emptyset$.
  \end{claim}
\begin{proof}
To prove this claim, fix $s\in p$ and $x\in G$.  Since $xC$ is
nowhere dense in $h\body{p_s}$ there exists $s_x\in p$ with $s\subseteq s_x$ 
and
$h\body{p_{s_x}}\cap xC=\emptyset$. By the continuity of the group operation
and the compactness of $xC$, 
there exists an open neighbourhood $U_x$ of $x$   with
$h\body{p_{s_x}}\cap U_xC=\emptyset$.  By  compactness,
finitely
many of these $U_x$ cover and so there is some  $N$  larger than
the length of the finitely many corresponding $s_x$ for each of
the finitely many $s$ in
$2^n \cap p$.    
\end{proof}

Now suppose  a sequence $(p^n,k_n)$ has been constructed satisfying
for every $n<\omega$:
\begin{enumerate}
 \item $p^n\in\Poset$, 
 \item $k_n<k_{n+1}<\omega$, 
 \item $p^{n+1}\leq p^n$,
 \item $2^{k_n}\cap p^{n+1}=2^{k_n}\cap p^{n}$, 
 \item any $s\in 2^{k_n}\cap p^{n+1}$ has incomparable extensions
 in $2^{k_{n+1}}\cap p^{n+1}$, and \label{splitting}
 \item for each $s\in p^{n+1}\cap 2^{k_n}$
and $x\in G$ there exists 
$t\in p^{n+1}\cap 2^{k_{n+1}}$ with $s\subseteq t$ and 
$h(\body{p^{n+1}_t})\cap xC=\emptyset$.  \label{nwdprop}
\end{enumerate}
Then $p=\bigcap_{n\in\omega}p^n$ is in $\Poset$.  Except for the last
condition \ref{nwdprop} 
 the fusion conditions  described are identical to those used by 
Baumgartner and Laver in \cite{MR81a:03050} for Sacks forcing.  The last
condition guarantees that for every $x\in G$ and $s\in p$ that
$xC$ is nowhere dense in $h(\body{p_s})$.  

In order to verify the properties (a), (b), and (c) note that
this fusion property is close enough to a property-A forcing to 
see easily
that $\Poset$ is proper.  Next it will be shown that forcing with $\Poset$
 is $\omega^\omega$-bounding.
Suppose 
$$p^0\forces {\tau:\omega\to Ordinals}.$$
Construct a fusion sequence $p^n,k_n$ as follows.  At stage
$n$ given $p^n,k_n$ find $p^{n+1}\leq p^n$ 
such that for every
$s\in 2^{k_n}\cap p^{n}$ the condition $p^{n+1}$ decides
$\tau(n)$.  Hence there exists a finite set $F_n$ such that
$$ p^{n+1}\forces {\tau(n)\in \check{F}_n}.$$
Now pick $k_{n+1}$ large enough to satisfy Condition \ref{splitting} and
Condition \ref{nwdprop}.  The fusion $p=\bigcap_{n\in\omega}p^n$ then
satisfies:
$$p\forces{\forall n\;\; \tau(n)\in \check{F}_n}.$$ 
Lastly, if $G$ is $\Poset_{(h,C,{\mathbb G})}$-generic over 
a ground model $V$ and $x\in 2^\omega$ is the real determined
by $G$ --- in other words, $x\in \body{p}$ for all $p\in G$ ---  then
$z=h(x)$ has the property that $z\notin XC$ where $X=V\cap G$.  This
follows from any easy density argument, since if $p\in\Poset$ and
$s\in p$ then $p_s\in\Poset$.

To prove the Theorem note that the countable support iteration of
proper forcings which are $\omega^\omega$-bounding is 
$\omega^\omega$-bounding (Shelah, see \cite{MR96k:03002} 6.3.5).
Hence a countable support $\omega_2$-iteration over a model of 
GCH will satisfy ${\mathfrak d}=\aleph_1$.
By dovetailing it is easy to arrange that every infinite compact
group and closed nowhere dense subset in the final model is forced
with (or rather a code for such) at unboundedly many stages in the iteration.
It follows by standard arguments that in the final model 
$\cov_{\mathbb G}=\aleph_2$ for every
infinite compact group ${\mathbb G}$.

\end{proof}

\begin{quest}
Is it consistent to have a compact group ${\mathbb G}$ such that 
$\cov_{\mathbb G}> \eq$?  
\end{quest}

\begin{quest}
Is it true that for any infinite compact group ${\mathbb G}$ that
$\cov_{\mathbb G}\geq \eq$?  
\end{quest}

\begin{quest}
Is it true that for every non-discrete Polish group ${\mathbb G}$
that $\cov_{\mathbb G}=\eq$ or $\cov_{\mathbb G}=\cov({\mathcal M})$?
\end{quest}

The question of general group actions has only been hinted at but
there are many questions concerning these as well. For example, it is
easy to see that if $\alpha_n$ is the natural action  of the isometry
 group  on $\Reals^n$ and $n\leq m$ then
$\cov_{\alpha_n} \leq \cov_{\alpha_m}$. However the following
question is unanswered.

\begin{quest}
  Is it true that $\cov_{\alpha_n} = \cov_{\alpha_m}$ for all
  $m$ and $n$?
\end{quest}


\begin{thebibliography}{1}

\bibitem{MR96k:03002}
Tomek Bartoszy{\'n}ski and Haim Judah.
\newblock {\em Set theory}.
\newblock On the structure of the real line.
\newblock A K Peters Ltd., Wellesley, MA, 1995.

\bibitem{MR81a:03050}
James~E. Baumgartner and Richard Laver.
\newblock Iterated perfect-set forcing.
\newblock {\em Ann. Math. Logic}, 17(3):271--288, 1979.

\bibitem{MR98d:54068}
Howard Becker and Alexander~S. Kechris.
\newblock {\em The descriptive set theory of {P}olish group actions}, volume
  232 of {\em London Mathematical Society Lecture Note Series}.
\newblock Cambridge University Press, Cambridge, 1996.

\bibitem{MR2000k:03097}
Greg Hjorth.
\newblock {\em Classification and orbit equivalence relations}, volume~75 of
  {\em Mathematical Surveys and Monographs}.
\newblock American Mathematical Society, Providence, RI, 2000.

\bibitem{LT}
Lindenstrauss, Joram; Tzafriri, Lior. Classical Banach spaces. I. Sequence
spaces. Ergebnisse der Mathematik und ihrer Grenzgebiete, Vol. 92.
Springer-Verlag, Berlin-New York, 1977. xiii+188 pp. ISBN: 3-540-08072-4

\bibitem{MR84e:03058a}
Arnold~W. Miller.
\newblock Some properties of measure and category.
\newblock {\em Trans. Amer. Math. Soc.}, 266(1):93--114, 1981.

\bibitem{MR84e:03058b}
Arnold~W. Miller.
\newblock Corrections and additions to: ``{S}ome properties of measure and
  category''.
\newblock {\em Trans. Amer. Math. Soc.}, 271(1):347--348, 1982.

\bibitem{vanmill}	
van Mill, J.; Infinite-dimensional topology. Prerequisites and introduction.
North-Holland Mathematical Library, 43. North-Holland Publishing Co.,
Amsterdam, 1989. xii+401 pp. ISBN: 0-444-87133-0.

\end{thebibliography}

\appendix

\newpage
\begin{center}
Appendix \\ 
\end{center}

This is not intended for publication but only for the electronic version.

\bigskip \noindent Details of the proof of Proposition \ref{torus} (2).
\bigskip

\begin{propo}\label{product}
  If ${\mathbb G}_n$ is a nontrivial finite group for each
  $n \in \Naturals$  and ${\mathbb G}=\prod_{n=0}^\infty{\mathbb G}_n$ 
  with the product topology then $\cov_{\mathbb{G}} = \eq$.
\end{propo}
\begin{proof}
To begin, note that if $f \leq g$ then $\eq(f)\geq \eq(g)$.
Choose $f:\Naturals \to  \Naturals$ such that $\eq(f) = \eq$. 
Now let $\{k_n\}_{n=0}^\infty$ be an increasing sequence of integers
such that $g(n) = \prod_{i=k_n}^{k_{n+1} - 1}|{G}_i| \geq f(n)$ for
each $n$. Clearly $\eq(g) = \eq$ as well.
Now let 
$C=\SetOf{x\in {G}}{(\forall n\in \Naturals)(\exists
i \in [k_n,k_{n+1}))\ x(i)\neq e_i}$ 
and note that  $C$ is closed nowhere dense.
Using the definition of $\eq(g)$ it is possible to find
$X\subseteq {G}$ such that for all $y\in {G}$ there is $x\in X$ such that 
$y\restriction [k_n,k_{n+1})\neq x \restriction [k_n,k_{n+1})$
for all $n$. Therefore $\prod_{n=0}^\infty (x(n) y(n)) \in C$; in
other words $XC = G$ and so $\cov_{\mathbb{G}} \leq \eq$.

On the other hand, suppose that
$C\subseteq G$ is closed nowhere dense 
and there is $W\subseteq G$ such that $|W| < \eq$ and $W C = G$.
Then it
is possible to find disjoint intervals of integers
$\{I_n\}_{n=0}^\infty$ such that, letting $\pi_n:{G} \to \prod_{i\in
  I_n}G_i$ be the projection map,  for all $n$ there is $x_n\in \prod_{i\in
  I_n}G_i$ such that $\pi_n^{-1}\{x_n\}\cap C = \emptyset$.
Let $$C^* = \SetOf{y\in \prod_{n=0}^\infty\prod_{i\in I_n}{G}_i}
                  {(\forall n)\ y(n)  \neq x_n}$$
and  let $g(n) = \prod_{i\in I_n}|{G}_i|$. 
Let $\{g_m\}_{m=1}^{g(n)}$ enumerate $\prod_{i\in I_n}{G_i}$. For
$w\in W$ let $f_w:\Naturals \to \Naturals$ be 
defined by $f_w(n) = j$ if and only if $ x_n  \pi_n(w) = g_j$.

Since $|W| < \eq(g)$ it is possible to find $h:\Naturals\to \Naturals$
such that for all $w\in W$
there is some $n \in \Naturals$ such that $f_w(n) = h(n)$.
Let $z = \prod_{n=0}^\infty g_{h(n)}$. It suffices to note that
$z\pi_n(w)^{-1}\notin C^*$ for all $w \in W$ because it then follows
that
$\pi^{-1}\{z\} \{w^{-1}\}\cap C = \emptyset$ for all $w\in W$
and so $W  C \neq G$. To see that $z \pi_n(w)^{-1}\notin
C^*$ let $n$ be such that $h(n) = f_w(n)$. Then $x_n  \pi_n(w)
 =
g_{h(n)} = z\restriction I_n$. Hence $ z \pi_n(w)^{-1}
= x_n$ and so $z  \pi_n(w)^{-1}\notin C^*$.
\end{proof}

\begin{propo}
For any integer $n$ the equality $\cov_{\mathbb{T}^n} = \eq$ holds.
\end{propo}

\begin{proof}
While there is a continuous
and open mapping from $2^\Naturals$ onto the circle $\mathbb T$, this mapping
is not a homomorphism.  Nevertheless, it is sufficiently close to a
homomorphism to be able to do the proof.

This follows the pattern of Proposition~\ref{product} taking care of  
carry digits. To be more precise, given $\psi:\Naturals\to            
\Naturals$ such that $\eq(\psi) = \eq$ let                      
$\{k_n\}_{n=0}^\infty$ be an increasing sequence of integers such     
that $\prod_{i=k_n}^{k_{n+1} - 1}|{G}_i| \geq 2\psi(n)$ for each      
$n$. Consider $\mathbb{T} $ to be $\Reals/\Integers$ so that          
$\mathbb{T}$ is identified with $[0,1]$.  For $t\in [0,1]$ let        
$x_t:\Naturals \to 2$ be chosen so that $\sum_{i=0}^\infty            
x_t(i)/2^i = t$ and there are infinitely many $i$ such that $x_t(i)   
= 0$.  Now let $$C=\SetOf{t\in [0,1]}{(\forall n\in                   
  \Naturals)(\exists i \in [k_n,k_{n+1}))\ x(i)\neq 0}$$              
and note                                                              
that $C$ is closed nowhere dense.  Now, for any two functions $f:     
[k_n,k_{n+1}-1) \to 2$ and $g: [k_n,k_{n+1}-1) \to 2$ define          
$f\equiv_n g$ if and only if there is some $j$ such that              
\begin{itemize}
 \item $f(i)\neq g(i)$ for $i \geq j$
 \item $f(i) = g(i)$ for $i < j$
 \item $f\restriction [j+1,k_{n+1})$ is constant ( and hence so is
   $g\restriction [j+1,k_{n+1})$).
\end{itemize}
It is immediate that $\equiv_n$ is an equivalence relation whose
equivalence classes are all pairs. Let $S_n$ be the set of equivalence
classes of $\equiv_n$. Hence $\phi(n) = |S_n| \geq \psi(n)$.

Using the definition of $\eq({\phi})$ it is possible to find
$Y\subseteq \prod_{n=0}^\infty S_n$ such that for all $t\in [0,1]$
there is $y\in Y$ such that $x_t\restriction [k_n,k_{n+1})\not\equiv_n
y(n)$ for all $n$. For each $y\in Y$ let $y^*$ be a choice function
which selects an element of $y(n)$ for each $n$ and let $\bar{y} =
\sum_{n=0}^\infty\sum_{i=k_n}{k_{n+1}-1}y^*(i)/2^i$. Observe that
$x_t\restriction [k_n,k_{n+1})\not\equiv_n x_{\bar{y}}\restriction
[k_n,k_{n+1})$ and hence, regardless of whether digits are carried
from the right, $x_{t - \bar{y}}\restriction [k_n,k_{n+1})$ is not
identically $0$. In other words, $t\in \bar{y} + C$.

On the other hand, suppose that $C\subseteq [0,1]$ is closed nowhere dense
and there is $W\subseteq [0,1]$ such that $|W| < \eq$ and $W + C =
[0,1]$.  Then it is possible to find disjoint intervals of integers
$\{I_n\}_{n=0}^\infty$ such that, letting $\equiv_n$ be the
equivalence relation defined on functions from $I_n$ to 2 and $S_n$ be
as before,
for all $n$ there is $[s_n]_{\equiv_n}\in S_n$ such that 
$$\SetOf{t\in [0,1]}{x_t\restriction I_n\equiv_n s_n}\cap C = \emptyset .$$
Let $$C^* = \SetOf{t\in[0,1]}{(\forall n)\ x_t\restriction I_n \not\equiv_n s_n}$$
and let $g(n) = 2^{|I_n|}$.  Let $\{g_m\}_{m=1}^{g(n)}$ enumerate
$2^{I_n}$. For $w\in W$ let $f_w:\Naturals \to
\Naturals$ be defined by $f_w(n) = j$ if and only if $  g_j - s_n =
x_w\restriction I_n $.

Since $|W| < \eq(g)$ it is possible to find $h:\Naturals\to \Naturals$
such that for all $w\in W$
there is some $n \in \Naturals$ such that $f_w(n) = h(n)$.
Now argue as in Proposition~\ref{product}.
\end{proof}

\bigskip

The following is another example of a corollary to \ref{c:g}.

\begin{corol}
  If $X$ is an orientable 2-manifold then
$\cov_{{\mathbb H}(X)} = \cov({\mathcal M})$.
\end{corol}
\begin{proof} 
An orientable 2-manifold is homeomorphic to the surface of the unit
  cube with holes drilled along axes perpendicular to the $xy$-plane.  
  Use Corollary~\ref{c:g} applied to the mapping which projects the
  surface to the $z$-axis.
\end{proof}

\end{document}